\documentclass{conm-p-l}

\newcommand{\<}{\kern.0833em}

\newtheorem{theorem}{Theorem}
\newtheorem{lemma}[theorem]{Lemma}
\newtheorem{corollary}[theorem]{Corollary}
\newtheorem{proposition}[theorem]{Proposition}

\newtheorem{question}[theorem]{Question}

\makeatletter
\newcommand{\xlabel}{\stepcounter{equation}
  \gdef\@currentlabel{\p@equation\theequation}{\rm(\@currentlabel)}}
\makeatother
\newenvironment{xlist}
  {\begin{list}{\xlabel}
    {\setlength{\rightmargin}{20pt}
     \setlength{\leftmargin}{37pt}
     \setlength{\labelsep}{20pt}
     \setlength{\labelwidth}{20pt}}}
  {\end{list}}
\makeatother

\begin{document}

\title[Two Statements about Infinite Products that Are Not Quite True]
{Two Statements about Infinite Products\\
that Are Not Quite True}
\keywords{infinite product module, inverse limit module, dichotomy of
finite or uncountable generation, homomorphism to infinite direct sum,
left perfect ring; infinite symmetric group, finite generation of
power or ultrapower of an algebra over the diagonal subalgebra.}

\author{George M. Bergman}
\dedicatory{To Don Passman, on his 65th birthday}
\address{Dept. of Mathematics\\
University of California\\
Berkeley, CA 94720-3840, USA}
\email{gbergman@math.berkeley.edu}

\subjclass[2000]{Primary: 08B25, 16D70.
Secondary: 03C20, 16P70, 16S50, 20B30, 20M20, 22A05}.
\date{}


\thanks{URLs of this preprint: 
http://math.berkeley.edu/{$\!\sim$}gbergman%
/papers/P\_\<vs\_cP.\{tex,dvi,ps,\linebreak[0]pdf\}~ (i.e.,
four versions, .../.tex, .../.dvi, .../.ps and .../.pdf).
\\
\indent
Comments, corrections, and related references are welcomed, as always!
}
\maketitle


\section{Introduction.}\label{S.intro}
The first half of this note concerns modules; so
let $R$ be a nonzero associative ring with unit.
A countably infinite direct product of copies of $R,$ generally
regarded as a left
$\!R\!$-module, will be denoted $R\<^\omega;$ the corresponding direct
sum, i.e., the free left $\!R\!$-module of countably infinite rank,
will be written $\bigoplus_\omega R.$

Here, now, are the two not-always-true statements of the title:

\begin{xlist}\item\label{x.cP!->P}
There is no surjective left $\!R\!$-module homomorphism
$\bigoplus_\omega R\rightarrow R\<^\omega.$
\end{xlist}
\begin{xlist}\item\label{x.P!->cP}
There is no surjective left $\!R\!$-module homomorphism
$R\<^\omega\to\bigoplus_\omega R.$
\end{xlist}

In \S\ref{S.cegs} we will note classes of rings $R$ for
which each of these statements fails.
In \S\S\ref{S.diag}-\ref{S.top}, however, we will see
that~(\ref{x.cP!->P}) holds, i.e., $R\<^\omega$ requires
uncountably many generators as a left $\!R\!$-module, unless
$R\<^\omega$ is {\em finitely} generated, and that~(\ref{x.P!->cP})
holds unless $R$ has descending chain condition on finitely
generated right ideals.

From the above assertion regarding~(\ref{x.cP!->P}),
and the statement of~(\ref{x.P!->cP}), it is
easy to see that for every $R,$ at least one
of~(\ref{x.cP!->P}),~(\ref{x.P!->cP}) holds.
We shall also see that the above restriction on
rings for which~(\ref{x.P!->cP}) fails implies that for every
$R,$ either~(\ref{x.P!->cP}) or the statement
\begin{xlist}\item\label{x.P!c>cP}
There is no embedding of left $\!R\!$-modules
$R\<^\omega \rightarrow \bigoplus_\omega R.$
\end{xlist}
holds.
(I did not count~(\ref{x.P!c>cP}) among the ``not quite true''
statements of the title, because I have no general results showing
that it is ``nearly'' true; i.e.,
that its failure implies strong restrictions on $R.)$

The result asserted above in connection with~(\ref{x.cP!->P})
will in fact be proved with $R\<^\omega$ replaced by $M^\omega$
for $M$ any $\!R\!$-module, while
the result on~(\ref{x.P!->cP}) will be obtained with $R\<^\omega$
generalized to the inverse limit of any countable inverse system
of finitely generated $\!R\!$-modules and surjective homomorphisms.

In~\S\ref{S.gen} we change gears:  We will note that our proof of the
result about~(\ref{x.cP!->P}) generalizes to the context of general
algebra (a.k.a `universal algebra'), and in \S\ref{S.Sym}, we
deduce from this that if $A$ is an algebra such
that $A^\omega$ is countably generated over the diagonal
image of $A,$ then it is finitely generated over that image.
We will then show that the monoid and group of all maps, respectively
all invertible maps, of an infinite set $\Omega$ into itself
have this finite generation property, and obtain results
on arbitrary algebras with this and related properties

I am grateful to P.\,M.\,Cohn, T.\,Y.\,Lam, B.\,Osofsky and
G.\,Sabbagh for pointing
me to some related literature, and to H.\,W.\,Lenstra,~Jr.\ and
T.\,Scanlon for contributions that will be noted below.

\section{Counterexamples.}\label{S.cegs}
It is easy to find rings for which statement~(\ref{x.P!->cP}) fails:
If $R$ is a division ring, then $R\<^\omega$ is
infinite-dimensional as a left $\!R\!$-vector-space, hence admits a
homomorphism onto $\bigoplus_\omega R.$
More generally, if $R$ is a quasi-Frobenius ring, then
the submodule $\bigoplus_\omega R \subset R\<^\omega,$ being free,
is injective~\cite[first paragraph]{CF}, so
$R\<^\omega$ can be retracted onto it, and again~(\ref{x.P!->cP}) fails.
The result to be proved in \S\ref{S.top}, that for~(\ref{x.P!->cP})
to fail $R$ must have descending chain condition on finitely
generated right ideals, shows that any example of failure
of~(\ref{x.P!->cP}) is fairly close to these cases.

Counterexamples to~(\ref{x.cP!->P}) are less evident.
To construct these, we begin by noting that given left modules $M$
and $N_i$ $(i$ ranging over an index set $I)$ over a ring $K,$ we have
\begin{xlist}\item\label{x.HomO+,}
${\rm Hom}_K(\bigoplus_I N_i,\,M)\,\cong\,\prod_I {\rm Hom}_K(N_i,\,M)$
as {\em right} $\!{\rm End}_K(M)\!$-modules
\end{xlist}
and
\begin{xlist}\item\label{x.Hom,*P}
${\rm Hom}_K(M,\ \prod_I N_i)\,\cong\,\prod_I {\rm Hom}_K(M,\,N_i)$
as {\em left} $\!{\rm End}_K(M)\!$-modules.
\end{xlist}
Indeed, the bijective correspondences follow from the universal
properties of the direct sum in~(\ref{x.HomO+,}) and the
direct product in~(\ref{x.Hom,*P});
it remains only to note that the left $\!K\!$-module
$M$ is a right $\!{\rm End}_K(M)\!$-module, and that this module
structure carries over to the hom-sets of~(\ref{x.HomO+,}), while it is
turned by contravariance into left module structures on the
hom-sets of~(\ref{x.Hom,*P}).
Here we have followed the standard convention of writing
homomorphisms of left modules on the right of their arguments,
and composing them accordingly (cf.~\cite{GMB.rl}), and we shall
continue to do so below.
We can now get our examples.

\begin{lemma}\label{L.R=R^*k}
Let $K$ be a ring and $\kappa$ a cardinal.
Then if $M$ is either a {\em right} $\!K\!$-module which satisfies
\begin{xlist}\item\label{x.M=(+)M}
$M~\cong~\bigoplus_\kappa M$
\end{xlist}
{\rm (}for instance, if $\kappa$ is infinite and $M$ a right
$\!K\!$-module of the form $\bigoplus_\kappa N),$
or a {\em left} $\!K\!$-module which satisfies
\begin{xlist}\item\label{x.M=M^}
$M~\cong~M^\kappa$
\end{xlist}
{\rm (}for instance, if $\kappa$ is infinite and $M$ a left
$\!K\!$-module of the form $N^\kappa),$
and if, in either case, we take $R~=~{\rm End}_K(M),$ then
\begin{xlist}\item\label{x.R=R^*k}
$R~\cong~R^\kappa$ as {\em left} $\!R\!$-modules.
\end{xlist}
\end{lemma}\begin{proof}
In~(\ref{x.HomO+,}) (with left and right interchanged),
respectively~(\ref{x.Hom,*P}) (as it stands), take $I=\kappa,$ put
the given module $M$ in the role of both the $M$ and all the~$N_i,$
and simplify the left hand side using~(\ref{x.M=(+)M}),
respectively~(\ref{x.M=M^}).
\end{proof}

A statement which embraces both cases of Lemma~\ref{L.R=R^*k}
is that if $\mathcal C$ is an $\!\mathcal{A}b\<\!$-cat\-e\-gory
with an object $M$ that is a $\!\kappa\!$-fold coproduct or
product of itself, then its endomorphism ring in $\mathcal C$ (or the
opposite of that ring, depending on which of these cases one is in,
and one's choice of how morphisms are composed in $\mathcal C)$
satisfies~(\ref{x.R=R^*k}).\vspace{6pt}

Clearly, when~(\ref{x.R=R^*k}) holds,~(\ref{x.cP!->P})
and~(\ref{x.P!c>cP}) fail.\vspace{6pt}

Let me say here how this subject came up.
H.\,W.\,Lenstra~Jr., in connection with a course he was teaching
in Leiden, e-mailed me asking whether, for a nontrivial ring $R,$
one could have
\begin{xlist}\item\label{x.hwl=?}
$R\<^\omega~\cong~\bigoplus_\omega R$
\end{xlist}
as left $\!R\!$-modules.
I replied that this was impossible because $R\<^\omega$ could
not be countably generated.
He pointed out the first case of Lemma~\ref{L.R=R^*k} (for $M$ an
abelian group, i.e., $K=\mathbb Z),$ which shows the contrary.
I then thought further and found the argument of the next section,
showing that~(\ref{x.hwl=?}) nonetheless cannot occur.

I have not been able to find in the literature any occurrence of
the idea of Lemma~\ref{L.R=R^*k} for infinite $\kappa,$
so it appears that the result belongs to Lenstra.
However, two less elementary examples of similar phenomena
were proven earlier.
J.\,D.\,O'Neill \cite{JO'Nrg} constructs for any $\kappa>1$ a ring $R$
without zero-divisors such that $R^\kappa\cong R^2$ as left modules,
and as noted in the MR review of that paper, W.\,Stephenson
showed in~\cite{WSt} that for any non-right-Ore ring $S$
without zero divisors, the right quotient ring $R$ of $S$
satisfies~(\ref{x.R=R^*k}) for any cardinal $\kappa$ such that
$S$ has $\geq\kappa$ right linearly independent elements.
That result is, in fact, an instance of the generalization
noted immediately after the proof of Lemma~\ref{L.R=R^*k} above.
For the right quotient ring of $S$ is the endomorphism
ring of $S$ in the $\!\mathcal{A}b\<\!$-category
$\mathcal C$ whose objects are right $\!S\!$-modules, but where
$\mathcal C(M,N)$ is the set of morphisms from essential
submodules of $M$ into $N,$ modulo the relation that
identifies morphisms which agree on essential submodules; and for
a non-right-Ore ring $S$ without zero-divisors having
$\geq\kappa$ right linearly independent elements,
the free right $\!S\!$-module on one generator has an essential
submodule consisting of a direct sum of $\geq\kappa$ copies of~$S,$
leading to an isomorphism~(\ref{x.M=(+)M}) in that category.

It is easy to give other sorts of examples where~(\ref{x.P!c>cP}) fails.
For instance, if $R$ is a noncommuting formal power series ring
on a $\!\kappa\!$-tuple of indeterminates over some ring $K,$ with no
restriction on the number of monomials allowed in each degree, then
its ideal of elements with constant term zero is isomorphic as a
left module to $R^\kappa,$ so $R^\kappa$ embeds in $R.$
On the other hand, examples where~(\ref{x.cP!->P}) fails are
harder to come by, so I will record a slight extension of the class of
examples given by Lemma~\ref{L.R=R^*k}.

First note that in any ring having a module
isomorphism~(\ref{x.R=R^*k}), e.g., a ring as in that lemma,
$R$ will contain elements $f_i$ $(i\in\kappa)$
such that {that} isomorphism is given by
\begin{xlist}\item\label{x.R->R^*k}
$r~\mapsto~(r f_i)_{i\in\kappa}.$
\end{xlist}
(E.g., in a case arising from a right module
isomorphism~(\ref{x.M=(+)M}), $(f_i)_{i\in\kappa}$
can be any family of one-to-one endomorphisms of
$M$ such that $M=\bigoplus_\kappa f_i(M).$
In general, the $f_i$ will be the components of the image
under~(\ref{x.R=R^*k}) of $1\in R.)$
But the surjectivity of~(\ref{x.R->R^*k}) is preserved
on replacing $R$ by any homomorphic image; so such a homomorphic
image will again be a ring for which~(\ref{x.cP!->P}) fails.
To get rings $R$ as in Lemma~\ref{L.R=R^*k} which have ideals $I$ such
that $R/I$ does not (so far as I know) itself satisfy~(\ref{x.R=R^*k}),
one can (a)~let $K$ be a field, $R$ the endomorphism ring of an
infinite-dimensional $\!K\!$-vector-space $V,$ and $I$ the ideal of
finite-rank endomorphisms of $V,$ or (b)~let $K$ be a ring having
a non-finitely-generated $\!2\!$-sided ideal $I_0,$ let $R=
{\rm End}_K(\bigoplus_\kappa K),$ and let $I$ be the ideal of $R$
generated by the diagonal image of $I_0.$
Because $I_0$ is non-finitely-generated, $I$ is generally smaller than
the kernel of the natural map ${\rm End}_K(\bigoplus_\kappa K)
\rightarrow {\rm End}_{K/I_0}(\bigoplus_\kappa K/I_0)$
(it does not contain endomorphisms whose components all lie in
$I_0,$ but do not all lie in any finitely generated subideal of $I_0),$
so $R/I$ is not simply the latter ring, which would just be
another example of the construction of Lemma~\ref{L.R=R^*k}.

It is interesting to note that a ring $R$ can
satisfy~(\ref{x.R=R^*k}) simultaneously as a right and as a left
$\!R\!$-module (though necessarily by different bijections).
Namely, if $K$ is a division ring and $\kappa$ an infinite
cardinal, then the left vector-space $M=K^\kappa$ not
only satisfies $M^\kappa\cong M,$ but also,
being a $\!K\!$-vector-space of dimension at least $\kappa,$ satisfies
$\bigoplus_\kappa M\cong M;$ hence combining Lemma~\ref{L.R=R^*k}
and its dual, we get the desired right and left module isomorphisms.

We now turn to results showing that despite these instances
of~(\ref{x.R=R^*k}), no nontrivial ring satisfies~(\ref{x.hwl=?}).

\section{A diagonal argument.}\label{S.diag}
In this section, we shall restrict ourselves, for simplicity of
presentation, to the case $\kappa=\omega.$
(In~\S\ref{S.gen}, in addition to passing from module theory
to general algebra, we will give the corresponding results
for a general infinite cardinal $\kappa.)$

The hypothesis of the next lemma may seem irrelevant to the
question at hand, but we shall see that in this case,
appearances are deceiving.
\begin{lemma}\label{L.Pinfg}
Let $R$ be a ring and $(M_i)_{i\in\omega}$ a family of
non-finitely-generated left $\!R\!$-modules.
Then $\prod_{i\in\omega}M_i$ is not countably generated.
\end{lemma}\begin{proof}
It will suffice show that for any countable family of elements
$x_j\in\prod_{i\in\omega} M_i$ $(j\in\omega),$ we can construct
an element $y$ not in the submodule generated by the $x_j.$
We do this by a diagonal construction:~ For each $i\in\omega,$
the assumption that $M_i$ is not finitely generated allows us
to take for the $\!M_i\!$-component of $y$ an element of $M_i$ not
in the span of the $\!M_i\!$-components of $x_0,\dots,x_i.$
If $y$ were in the span of all the $x$'s, it would be in the
span of finitely many of them, say $x_0,\dots,x_i.$
But looking at its $\!M_i\!$-component, we get a contradiction.
\end{proof}

The relevance of that lemma can now be seen in the proof of

\begin{theorem}\label{T.fgorunc}
Let $M$ be a left module over a ring $R.$
Then the left $\!R\!$-module $M^\omega$ is either finitely generated,
or not countably generated.
\end{theorem}\begin{proof}
If $M^\omega$ is not finitely generated, then the above
lemma shows that $(M^\omega)^\omega$ is not countably generated.
But $(M^\omega)^\omega\cong M^{\omega\times \omega}\cong M^\omega,$
so $M^\omega$ is non-countably-generated, as claimed.
\end{proof}

So Lenstra's question is answered:

\begin{corollary}\label{C.hwl_ans1}
No nonzero ring $R$ can satisfy $R\<^\omega\cong\,\bigoplus_\omega R$
as left $\!R\!$-modules.
\end{corollary}\begin{proof}
$\bigoplus_\omega R$ is neither finitely generated nor
non-countably-generated, while we have seen that $R\<^\omega$
must have one of these properties.
\end{proof}

Let us denote the cardinality of a set $X$ by $|X|.$
I do not know the answer to
\begin{question}\label{Q.c'uum}
In the context of Theorem~\ref{T.fgorunc}, if $M^\omega$
is not finitely generated, must the least cardinality of a
generating set be $\geq 2^{\aleph_0}$?
\end{question}

If $|R\<|<|M|^{\aleph_0},$ in particular if $|R\<|<2^{\aleph_0},$
the answer is affirmative even without an explicit assumption
that $M^\omega$ is not finitely generated, as long as $M\neq 0;$
indeed, a generating set must have cardinality
$|M|^{\aleph_0},$ since a module of infinite cardinality cannot be
generated by a subset of smaller cardinality
over a ring of smaller cardinality.
Also, without any restriction on $|R\<|,$ the module $M^\omega$
will {\em contain} a direct sum of $2^{\aleph_0}$ copies of $M.$
This can be seen either by applying
\cite[Corollary~3.5]{ZB} to the discrete topology
on $\omega,$ or, alternatively, as follows:
Note that $M^\omega\cong M^\mathbb Q$ where $\mathbb Q$ is the set
of rational numbers; associate to every real number
$r$ the map $f_r: M\rightarrow M^\mathbb Q$ carrying each
$x\in M$ to the function on $\mathbb Q$ which has the
value $x$ at all rational numbers $q<r$ and $0$ at all $q\geq r,$
and verify that the sum of the resulting images of $M$ is direct.

For $R$ a division ring, this shows that
$R\<^\omega$ has dimension at least $2^{\aleph_0};$ but in this case
the Erd\H{o}s-Kaplansky Theorem \cite[Theorem~IX.2, p.246]{NJ} gives
the more precise result ${\rm dim}_R\,R\<^\omega=|R\<|^{\aleph_0}.$
This suggests that perhaps in
general the least cardinality of a generating set
of a module $M^\omega,$ if infinite,
must equal $|M|^{\aleph_0};$
but the following example shows that this is not so.
Let $R=R'\times R'',$ where $R'$ is a ring as in Lemma~\ref{L.R=R^*k}
and having cardinality greater
than $2^{\aleph_0},$ while $R''$ is any nontrivial
ring of cardinality $\leq 2^{\aleph_0}$ (e.g., a finite
ring) such that $R''^\omega$
is not finitely generated as a left $\!R\!$-module.
Then $R\<^\omega\cong{R'}^\omega\times{R''}^\omega$ can be generated
by a cyclic generator of ${R'}^\omega$ together with the $2^{\aleph_0}$
elements of ${R''}^\omega,$ hence by fewer than $|R|$ elements.
So though Question~\ref{Q.c'uum} is open, there is no obvious
stronger conjecture to make.

The argument by which Theorem~\ref{T.fgorunc} was obtained
from Lemma~\ref{L.Pinfg} can be applied to show uncountable
generation of other sorts of product modules.
For instance, if $(M_i)_{i\in\omega}$ is a family of finitely generated
modules such that the finite numbers of generators they require
is unbounded,
it is clear that their direct product cannot be finitely generated.
But such a family of modules can be partitioned into countably many
infinite subfamilies each having the same unboundedness property;
hence by Lemma~\ref{L.Pinfg}, its direct product is
in fact non-countably-generated.

This leaves open the case where the number of generators of
the $M_i$ is bounded.
By passing to matrix rings, one can reduce this to the cyclic case,
so we ask
\begin{question}\label{Q.bddgen}
Let $R$ be a ring and $(M_i)_{i\in\omega}$ a family of cyclic
left $\!R\!$-modules.
Must the $\!R\!$-module $\prod_{i\in\omega} M_i$ either be
finitely generated or require uncountably many generators?
\end{question}

If $(M_i)_{i\in\omega}$ is any family of cyclic
$\!R\!$-modules, let $V=\{S\,{\subseteq}\,\omega\,|\,\prod_{i\in S}
M^i$ is not finitely generated$\!\}.$
Whenever $V$ contains a union $S\cup S'$ of two sets,
it must clearly contain $S$ or $S',$
and if it contains infinitely many disjoint sets, an argument
like that just noted shows that $\prod_{i\in\omega} M_i$ is
uncountably generated.
From this it is not hard to show that if $(M_i)_{i\in\omega}$ is
a counterexample to Question~\ref{Q.bddgen}, $V$ must be a union
of finitely many nonprincipal ultrafilters on $\omega.$
Taking $S\subseteq\omega$ which belongs to exactly one of these
ultrafilters, and reindexing by $\omega,$ we get such a
counterexample where $V$ itself is an ultrafilter.
But I don't see how to go anywhere from there --
the obvious thought is ``look at the ultraproduct module
$(\prod_{i\in\omega} M_i)/V$\!'', but I see no reason why that
module would have to be non-finitely-generated,
let alone uncountably generated.

\section{A touch of topology, and some chain conditions.}\label{S.top}
When I sent him Theorem~\ref{T.fgorunc}, Lenstra noted that this
implied not only that the left $\!R\!$-modules $R\<^\omega$
and $\bigoplus_\omega R$ were nonisomorphic, but that
one could not map each surjectively to the other; i.e.,
in our present notation, that~(\ref{x.cP!->P})
and~(\ref{x.P!->cP}) could not fail simultaneously.
He then raised the question of whether for some $R$
the module $R\<^\omega$ could admit both injective and surjective
maps to $\bigoplus_\omega R;$ i.e., whether~(\ref{x.P!->cP})
and~(\ref{x.P!c>cP}) could both fail.

This is also impossible.
I know now that this fact can be obtained without too much work
from a result of S.\,U.\,Chase on homomorphisms from direct product
modules to direct sums~\cite[Theorem~1.2]{SUC3}.
However, I will prove some stronger statements, including a
generalization of Chase's result, Theorem~\ref{T.MNB-Chase} below
-- essentially, the convex hull of his theorem and the similar
result that I had obtained before learning of~\cite{SUC3}.

Before leaping into the proof, let us note that a general tool
in proving restrictions on homomorphisms $f$ out of an infinite
product module $M=\prod_i M_i$
is to assume those restrictions fail, and
construct an element $x\in M$ by specifying its values
on successive coordinates, in such a way that the properties of $f(x)$
would lead to a contradiction.
To ``control'' the effects of these successive specifications, one
generally takes the $\!i\!$th coordinate to lie in an additive subgroup
$I_i\,M_i\subseteq M_i,$ using smaller and smaller right ideals
$I_i\subseteq R,$ and calling (explicitly or implicitly) on the fact
that if a right ideal $I$ is finitely generated, and we modify an
arbitrary set of coordinates $x_i$ of $x$ by elements of $I\<M_i,$
then this modifies $f(x)$ by an element of the subgroup $I\<f(M).$

Both Chase's proof in~\cite{SUC3} (and his proofs of similar
results in \cite{SUC1} and \cite{SUC2}) and my original argument
used this method of successive approximation.
However, after a series of generalizations and reformulations,
in which the property of the direct product being used
was translated into a completeness condition with respect to an inverse
limit topology, I realized that {that} step
was essentially a repetition of the proof of the
Baire Category Theorem, and could be avoided by calling on that theorem.
Here is that part of the argument:
\begin{lemma}[{cf.\ \cite[Lemma~3.3.3]{EK}}]\label{L.Baire}
Let $G$ be a complete metrizable topological group,
and $(B_i)_{i\in\omega}$ a countable family of subgroups of $G$
such that $\bigcup_{i\in\omega} B_i=G.$
Then for some $i,$ the closure of $B_i$ in $G$ is open.
\end{lemma}\begin{proof}
The closed subgroups ${\rm cl}(B_i)$ again have
union $G,$ hence by the Baire Category Theorem, some
${\rm cl}(B_i)$ has nonempty interior, i.e., contains a
neighborhood in $G$ of one of its points $x.$
By translation, it contains a neighborhood of each
of its points, hence it is open.
\end{proof}

Although I spoke above of product modules, the same methods are
applicable, more generally, to inverse limits of countable systems
of modules, and we will prove our results for these.
Note that if a module $M$ is the inverse limit of a system
\begin{xlist}\item\label{x.invsys}
$\dots~\rightarrow M_i\rightarrow\ \dots\ \rightarrow M_2\rightarrow
M_1\rightarrow M_0$
\end{xlist}
with surjective connecting homomorphisms, then $M$ maps surjectively
to each $M_i.$
(The corresponding statement is not true of inverse limits over
uncountable partially ordered sets \cite[Example~10.4]{EK}!)
Hence each $M_i$ may be written $M/N_i,$ where the kernels form a chain
\begin{xlist}\item\label{x.Ni}
$N_0\,\supseteq N_1\,\supseteq\,\dots\,\supseteq\,N_i\,\supseteq\,
\dots\,,$ such that $\bigcap_{i\in\omega} N_i=\{0\},$ and $M$ is
complete in the $\!(N_i)\!$-adic topology.
\end{xlist}

(A countable {\em direct product} module $M=\prod_{i\in\omega} L_i$
is the case of~(\ref{x.invsys}) and~(\ref{x.Ni}) in which for all $i,$
$M_i=\prod_{j=0} ^{i-1}\,L_j$ and $N_i=\prod_{j\geq i}\,L_j.)$

The finitely generated right ideals $I_i$ forming the other ingredient
of the technique sketched above will come into the
picture through the following curious result.

\begin{lemma}[{\rm cf.\ \cite[Lemma~3.3.4]{EK}}\textbf{}]\label{L.IN}
Let $R$ be a ring, and $M$ a left $\!R\!$-module complete with
respect to a chain of submodules~{\rm (\ref{x.Ni})}.
Then for any chain of finitely generated right ideals of $R,$
\begin{xlist}\item\label{x.Ii}
$I_0~\supseteq~ I_1~ \supseteq~\dots~\supseteq~ I_i~ \supseteq~\dots\,,$
\end{xlist}
$M$ is also complete with respect to the chain of additive subgroups
\begin{xlist}\item\label{x.IN}
$I_0 N_0~\supseteq~I_1 N_1~\supseteq~\dots~\supseteq~
I_i N_i~\supseteq~\dots~.$
\end{xlist}
\end{lemma}\begin{proof}
Let $(x_i)$ be a sequence of elements of $M$ satisfying
\begin{xlist}\item\label{x.x+IN}
$x_{i+1}~\in~x_i+I_i\,N_i$ for all $i\in\omega.$
\end{xlist}
Since $I_i N_i\subseteq N_i,$ the sequence $(x_i)$
is Cauchy with respect to the $\!(N_i)\!$-adic topo\-logy, and so
has a limit $x$ in that topology.
Our desired conclusion will clearly follow if we can show that $(x_i)$
converges to $x$ in the $\!(I_i N_i)\!$-adic topology as well.
To do this it will suffice to show that for each $i_0\in\omega,$
\begin{xlist}\item\label{x.xi0+}
$x~\in~x_{i_0} + I_{i_0} N_{i_0}.$
\end{xlist}
Let us fix $i_0$ for the remainder of the proof, and
establish~(\ref{x.xi0+}).

We first note that~(\ref{x.x+IN}) entails the weaker statement
gotten by ignoring cases before the $\!i_0\!$th, and
replacing all the ideals $I_i$ $(i\geq i_0)$ with the
larger ideal $I_{i_0}:$
\begin{xlist}\item\label{x.x+Ii0N}
$x_{i+1}~\in~x_i+I_{i_0}\,N_i$ for $i\geq i_0.$
\end{xlist}
By assumption, $I_{i_0}$ has a finite generating set $S$ as a right
ideal, so~(\ref{x.x+Ii0N}) shows that for each $i\geq i_0$ we can write
\begin{xlist}\item\label{x.xi+1}
$x_{i+1}~=\ x_i+\sum_{s\in S}\,s\,y_{si},$ with $y_{si}\in N_i$
for each $s\in S.$
\end{xlist}
By completeness of $M$ in the $\!(N_i)\!$-adic topology,
for each $s\in S$ the
series $\sum_{i\geq i_0} y_{si}$ converges in that topology to
an element $y_s\in N_{i_0},$ so for each $s,$
$\sum_{i\geq i_0}\,s\,y_{si}$ converges in the
$\!(I_{i_0}N_i)\!$-adic topology to $s\,y_s\in I_{i_0} N_{i_0}.$
Summing over $S$ and adding to $x_{i_0},$ we conclude
from~(\ref{x.xi+1}) that the sequence $(x_i),$ converges in
that topology to
\begin{xlist}\item\label{x.xi0+Sum}
$x_{i_0}+\sum_S s\,y_s.$
\end{xlist}

But by assumption, the sequence $(x_i)$ converges in the
$\!(N_i)\!$-adic topology to $x,$ hence the limit~(\ref{x.xi0+Sum})
of that sequence in the stronger $\!(I_{i_0}N_i)\!$-adic must also be
$x,$ proving~(\ref{x.xi0+}), as required.
\end{proof}

Now consider a situation where we have a homomorphism $f$ from
an inverse limit module $M$ as above (for instance, $R\<^\omega)$ onto
the free module of countable rank, $\bigoplus_\omega R.$
For each $j\geq 0,$ the elements $x\in M$ such that $f(x)$
has no nonzero components after the $\!j\!$th component form
a submodule $B_j\subseteq M,$ and since there are no elements $x$
such that $f(x)$ has infinitely many nonzero components,
$\bigcup_{j\in\omega} B_j = M.$
In the next result we play the existence
of such a chain of submodules off against completeness.

By a {\em downward directed system $\mathcal F$}
of right ideals of $R,$ we shall mean
a set $\mathcal F$ of right ideals such that every pair of members
of $\mathcal F$ has a common lower bound in $\mathcal F.$

\begin{theorem}\label{T.MNB}
Let $R$ be a ring and $M$ a left $\!R\!$-module which has
a chain of submodules~{\rm (\ref{x.Ni})}; equivalently,
which is the inverse limit of a system~{\rm (\ref{x.invsys})} of
modules and surjective homomorphisms.
Suppose we are also given an ascending chain of submodules
\begin{xlist}\item\label{x.B_j}
$B_0~\subseteq~B_1~\subseteq\,\dots$ with union $M,$
\end{xlist}
and a downward directed system $\mathcal F$ of right ideals of $R.$
Then there exists $j^*\in\omega$ such that, writing
$q_{j^*}$ for the canonical map $M\rightarrow M/B_{j^*},$ we have
\begin{xlist}\item\label{x.IN/B}
The set $\{I\,q_{j^*}(N_k)\ |\ I\in\mathcal F,\ k\in\omega\}$
of additive subgroups of $M/B_{j^*}$ has a least member.
\end{xlist}
I.e., there exists $I^*\in\mathcal F$ and $k^*\in\omega$ such that
$I^*\,q_{j^*}(N_{k^*})\subseteq I\,q_{j^*}(N_k)$ for all
$I\in\mathcal F,\ k\in\omega.$
\end{theorem}\begin{proof}
Suppose~(\ref{x.IN/B}) is false.
We will begin by constructing recursively a sequence
$I_0\supseteq I_1\supseteq\dots$ of ideals from $\mathcal F,$ starting
with an arbitrary $I_0\in\mathcal F,$ and a
sequence $k_0<k_1<\dots$ of natural numbers, starting with $k_0=0.$
Assuming $I_0,\dots,I_{i-1}$ and $k_0,\dots,k_{i-1}$ constructed,
the additive subgroup
$I_{i-1}\,q_i(N_{k_{i-1}})\subseteq M/B_i$ is by assumption
{\em not} least among subgroups $I\,q_i(N_k)$
$(I\in\mathcal F,\,k\in\omega),$ so we can pick $I_i\subseteq I_{i-1}$
and $k_i>k_{i-1}$ giving a proper inclusion
\begin{xlist}\item\label{x.smaller}
$I_i\,q_i(N_{k_i})~\subset\ I_{i-1}\,q_i(N_{k_{i-1}}).$
\end{xlist}

Once these choices have been made for all $i,$ Lemma~\ref{L.IN} tells us
that $M$ is complete in the $\!(I_i N_{k_i})\!$-adic topology, hence by
Lemma~\ref{L.Baire}, there is some $j\in\omega$ such that the closure of
$B_j$ in that topology is also open.
Thus for some $i_0\in\omega,$ that closure contains
$I_{i_0} N_{k_{i_0}},$ which is easily seen to imply that
\begin{xlist}\item\label{x.allIN}
$B_j+I_i\,N_{k_i}~\supseteq\ I_{i_0} N_{k_{i_0}}$ for all $i\in\omega.$
\end{xlist}
This means that when $i\geq i_0,$ taking larger values of $i$
(and hence smaller subgroups $I_i\,N_{k_i})$
does not decrease the left-hand side of~(\ref{x.allIN}).
Also, note that~(\ref{x.allIN}), and hence
that conclusion, is preserved on increasing $j.$
We easily deduce
\begin{xlist}\item[]
$B_{j'} + I_i\,N_k~=\ B_{j'} + I_{i'}\,N_{k'}$
for all $i,i'\geq i_0,\ j'\geq j,\ k,k'\geq k_{i_0}.$
\end{xlist}
or as an equation in $M/B_{j'},$
\begin{xlist}\item[]
$I_i\,q_{j'}(N_k)~=\ I_{i'}\,q_{j'}(N_{k'})$
for all $i,i'\geq i_0,\ j'\geq j,\ k,k'\geq k_{i_0}.$
\end{xlist}
But the strict inclusion~(\ref{x.smaller}), for any
$i\geq\mathrm{max}(j,i_0{+}1),$ contradicts this equality,
completing the proof.
\end{proof}

The conclusion of the above theorem is rather complicated, with
quantification over
$\mathcal F,$ $(B_j)_{j\in\omega}$ and $(N_k)_{k\in\omega}.$
One can get a statement involving only the first two of these
if one assumes that $M$
is an inverse limit of {\em finitely generated} $\!R\!$-modules $M_i.$
\begin{corollary}\label{C.fg}
Suppose in the situation of Theorem~\ref{T.MNB} that the $\!R\!$-modules
$M_i=M/N_i$ of~{\rm (\ref{x.invsys})} are finitely generated.
Then there exists $j^*\in\omega$ such that
\begin{xlist}\item\label{x.IM/B}
the set of additive subgroups $\{I\,(M/B_{j^*})\ |\ I\in\mathcal F\}$
of $M/B_{j^*}$ has a least member.
\end{xlist}
{\rm (}This will also be the least member of the larger family
of additive subgroups described as in~{\rm (\ref{x.IN/B}).)}
\end{corollary}\begin{proof}
Let $j^*,$ $I^*,$ $k^*$ be as in Theorem~\ref{T.MNB}, and
note that the property asserted by that theorem is preserved under
increasing $j^*$ while leaving $I^*$ and $k^*$ unchanged.
Let $M/N_{k^*}$ be generated by the image
of a finite set $T\subseteq M.$
Since $\{B_j\}$ has union $M,$ we can assume
by increasing $j^*$ if necessary that $T\subseteq B_{j^*}.$
This says that the image of $B_{j^*}$ in $M/N_{k^*}$ is all of
$M/N_{k^*},$ equivalently, that $M=B_{j^*} + N_{k^*},$ equivalently,
that the image of $N_{k^*}$ in $M/B_{j^*}$ is all of $M/B_{j^*}.$
Substituting this into the final assertion of Theorem~\ref{T.MNB}
gives the desired conclusion.
\end{proof}

We can now prove our earlier
assertion about maps onto $\bigoplus_\omega R.$

\begin{corollary}\label{C.DCC}
Suppose $R$ is a ring such that the free left $\!R\!$-module
$\bigoplus_\omega R$ can be written as a homomorphic image of the
inverse limit $M$ of an inversely directed system~{\rm (\ref{x.invsys})}
of finitely generated left $\!R\!$-modules and surjective homomorphisms.
Then $R$ is left perfect, i.e., has descending chain condition
on finitely generated right ideals.

In particular, any ring $R$ for which
there exists a surjective left $\!R\!$-module
homomorphism $R\<^\omega\rightarrow\bigoplus_\omega R,$
i.e., for which {\rm (\ref{x.P!->cP})} fails, is left perfect.
\end{corollary}\begin{proof}
Given $M$ as above
and a surjective homomorphism $f:M\rightarrow \bigoplus_\omega R,$
let $B_i=f^{-1}(\bigoplus_0 ^{i-1}\,R)$ $(i\in\omega).$
Then $M=\bigcup B_i;$ moreover, each factor module $M/B_i$ is isomorphic
to $\bigoplus_{j\geq i} R,$ again a free left module of countable
(hence nonzero) rank,
so for every $i\in\omega,$ distinct right ideals
$I$ yield distinct abelian groups $I(M/B_i).$
The conclusion of Corollary~\ref{C.fg} therefore tells us
that every downward directed system of finitely generated right ideals
of $R$ has a least member.
Applied to chains of ideals, this says that $R$ has
DCC on finitely generated right ideals, as claimed.
The final assertion is clear.
\end{proof}

We noted at the beginning of \S\ref{S.cegs} that~(\ref{x.P!->cP})
always failed when $R$ was a quasi-Frobenius ring;
so the class of rings for which it fails lies between the
classes of quasi-Frobenius and left perfect rings.
It would be of interest to know whether it equals one of
those classes, and if not, to characterize it better.

H.\,Lenzing \cite[Proposition~2]{Lenzing} says
(when restated in left-right dual form) that if
$\bigoplus_{i\in\omega} R\subseteq R^\omega$
is a direct summand as left modules, then $R$ is semiprimary
with ascending chain condition on left annihilator ideals.
(Whether the converse holds is left open
\cite[sentence beginning at bottom of p.687]{Lenzing}.)
The failure of our condition~(\ref{x.P!->cP}) is equivalent to
the statement that $\bigoplus_{i\in\omega} R$ can be embedded in
{\em some} way as a direct summand in
$R^\omega;$ so it is natural to ask,

\begin{question}\label{Q.summand}
If $\bigoplus_{i\in\omega} R$ has some embedding as a direct summand
in $R^\omega,$ will the canonical copy
of $\bigoplus_{i\in\omega} R$ in $R^\omega$ be a direct summand?
\end{question}

Getting back to our original
goal, we can now obtain the result Lenstra asked about.

\begin{corollary}\label{C.hwl_ans2}
For a nonzero ring $R,$ statements~{\rm (\ref{x.P!->cP})}
and~{\rm (\ref{x.P!c>cP})} cannot fail simultaneously;
i.e., $R\<^\omega$ cannot be both mappable onto and
embeddable in $\bigoplus_\omega R.$
\end{corollary}\begin{proof}
By the observations following Question~\ref{Q.c'uum}, a power module
$R\<^\omega$ always contains a direct sum of $2^{\aleph_0}$ copies
of $R,$ hence if $R\<^\omega$ is embeddable in $\bigoplus_\omega R,$
that module will contain
a set $X$ of $2^{\aleph_0}$ left linearly independent elements.
Writing $\bigoplus_\omega R$ as  the union of the countable
ascending chain of submodules $\bigoplus_0 ^{i-1} R,$ we conclude that
for some $i,$ this submodule will contain $2^{\aleph_0}$ members
of $X;$ so in particular it contains $i\,{+}\,1$ members
of $X;$ in other words, a free left $\!R\!$-module
of rank $i$ contains a free submodule of rank $i\,{+}\,1.$
This leads to a left regular $i{+}1\times i$ matrix over $R,$ hence
to an $i{+}1\times i{+}1$ matrix with left linearly
independent rows, but with a zero column, hence to an element $A$
of the full $i{+}1\times i{+}1$ matrix ring over $R$ which is left
but not right regular in that ring.

On the other hand, if $R\<^\omega$ is mappable onto
$\bigoplus_\omega R,$ then by the preceding corollary, $R$ is left
perfect, hence so is its $i{+}1\times i{+}1$ matrix ring.
A generalization by Lam \cite[Exercise~21.23]{Lam} of a result of
Asano's says that the left regular elements, the right regular elements
and the invertible elements in such a ring coincide, so an element $A$
cannot be left but not right regular.

Hence the two conditions on $R$ are incompatible.
\end{proof}

Returning to Theorem~\ref{T.MNB}, let me now obtain from it a version
with the conclusion in a form closer to Chase's formulation.
Note that the hypothesis below is as in Theorem~\ref{T.MNB}, except that
instead of a countable ascending chain of submodules
$(B_i)_{i\in\omega},$ we assume given a
homomorphism from $M$ into a module $\bigoplus_{\alpha\in J} C_\alpha,$
with no restriction on the cardinality of the index set~$J.$

\begin{theorem}
[{\rm after Chase \cite[Theorem~1.2]{SUC3}}\textbf{}]\label{T.MNB-Chase}
Let $R$ be a ring and $M$ a left $\!R\!$-module which has
a descending chain of submodules~{\rm (\ref{x.Ni})}; equivalently,
which is the inverse limit of a system~{\rm (\ref{x.invsys})} of
modules and surjective homomorphisms, and let $\mathcal F$ be
a downward directed system of right ideals of $R.$
Suppose we are given a family $(C_\alpha)_{\alpha\in J}$ of left
$\!R\!$-modules, and
a homomorphism $f: M\rightarrow \bigoplus_{\alpha\in J}\,C_\alpha.$
For each $\beta\in J$ let $\pi_\beta: \bigoplus_{\alpha\in J}\,C_\alpha
\rightarrow C_\beta$ denote the $\!\beta\!$th projection map.

Then there exist $k^*\in\omega,$ a finite subset $J_0\subseteq J,$
and an $I^*\in\mathcal F$ such that
\begin{xlist}\item\label{x.pfN=cap}
$I^*\,\pi_\beta(f(N_{k^*}))~\subseteq\ \bigcap_{I\in\mathcal F}\,
I\,C_\beta$~ for all $\beta\in J-J_0.$
\end{xlist}
\end{theorem}\begin{proof}
Assuming the contrary, let us construct recursively
a chain of ideals $I_0\supseteq I_1\supseteq\dots$ in $\mathcal F$ and
a sequence of indices $\alpha_0,\,\alpha_1,\dots\in J,$ starting
with arbitrary $I_0\in\mathcal F$ and $\alpha_0\in J.$

Say $I_0\supseteq\dots\supseteq I_{i-1}$
and $\alpha_0,\dots,\alpha_{i-1}$ have been chosen for some $i>0.$
By assumption, the choices
$I^*=I_{i-1},\ k^*=i,\ J_0=\{\alpha_0,\dots,\alpha_{i-1}\}$ do
not satisfy~(\ref{x.pfN=cap}), hence we can choose some $\alpha_i\not\in
\{\alpha_0,\dots,\alpha_{i-1}\}$
such that $I_{i-1} \,\pi_{\alpha_i}(f(N_i))\not\subseteq
\bigcap_{I\in\mathcal F} I\,C_{\alpha_i}.$
This in turn says that we can choose some $I_i\in\mathcal F$ such that
\begin{xlist}\item\label{x.nonstop}
$I_{i-1}\,\pi_{\alpha_i}(f(N_i))~\not\subseteq\ I_i\,C_{\alpha_i}.$
\end{xlist}
Since this property is preserved under replacing $I_i$ by
a smaller ideal, and since $\mathcal F$ is downward directed,
we may assume that $I_i\subseteq I_{i-1}.$

After constructing these ideals and indices for all $i,$
we define for each $j\in\omega$
\begin{xlist}\item[]
$B_j~=\ \{x\in M\ |\ \forall i\geq j,
\ \pi_{\alpha_i}(f(x))\in\bigcap_{I\in\mathcal F} I\,C_{\alpha_i}\}.$
\end{xlist}
Note that for any $x\in M,$ the finite support of $f(x)$ contains
$\alpha_i$ for only finitely many $i,$ so taking $j$ greater
than the largest of these values, we see that $x\in B_j.$
Thus, the ascending chain of submodules $B_j$ has union $M.$

Applying Theorem~\ref{T.MNB}, with $\{I_i\ |\ i\in\omega\}$ in place
of $\mathcal F,$ we get $i^*,\,j^*,\,k^*$ such that in $q_{j^*}(M)=
M/B_{j^*},$
\begin{xlist}\item[]
$I_i(q_{j^*}(N_k))~=
\ I_{i^*}(q_{j^*}(N_{k^*}))$ for all $i\geq i^*,\ k\geq k^*.$
\end{xlist}
From the definition of $B_{j^*}$ and $q_{j^*},$ this implies that
\begin{xlist}\item[]
$I_i\,\pi_{\alpha_j}(f(N_k))\,+\,
\bigcap_{I\in\mathcal F} I\,C_{\alpha_i}
~=\ I_{i^*}\,\pi_{\alpha_j}(f(N_{k^*}))\,+\,
\bigcap_{I\in\mathcal F} I\,C_{\alpha_i}$
for all $i\geq i^*,\ j\geq j^*,\ k\geq k^*.$
\end{xlist}
Hence if we fix $j\geq j^*$ and vary
$i\geq i^*$ and $k\geq k^*,$ the left hand side above does not depend
on the latter two values.
It follows that if we take $i={\rm max}(i^*{+}1,j^*,k^*),$ we have
\begin{xlist}\item[]
$I_{i-1}\,\pi_{\alpha_i}(f(N_i))+
\bigcap_{I\in\mathcal F} I\,C_{\alpha_i}
~=\ I_i\,\pi_{\alpha_i}(f(N_i))+
\bigcap_{I\in\mathcal F} I\,C_{\alpha_i,}$
\end{xlist}
so
\begin{xlist}\item[]
$I_{i-1}\,\pi_{\alpha_i}(f(N_i))
\ \subseteq\ I_i\,\pi_{\alpha_i}(f(N_i))+
\bigcap_{I\in\mathcal F} I\,C_{\alpha_i}\\
\indent\subseteq\ I_i\,C_{\alpha_i}+\bigcap_{I\in\mathcal F}
I\,C_{\alpha_i}\ =~I_i\,C_{\alpha_i},$
\end{xlist}
contradicting~(\ref{x.nonstop}) and completing the proof.
\end{proof}

As with Theorem~\ref{T.MNB}, if one assumes the modules $M_i=
M/N_i$ finitely generated, one gets a simplified conclusion, with
$I^*\,\pi_\beta(f(M))$ in place of $I^*\,\pi_\beta(f(N_{k^*})).$
\vspace{6pt}

Is my development above an improvement on Chase's
proof of \cite[Theorem~1.2]{SUC3}?
It improves the result by
generalizing direct products to inverse limits, and principal
right ideals to finitely generated right ideals, but these changes
could have been made without significantly altering his argument.
The above development, with the auxiliary lemmas, and the derivation of
Theorem~\ref{T.MNB-Chase} via Theorem~\ref{T.MNB}, is
longer than Chase's.
The best I can say is that it provides alternative perspectives on
what underlies these results, complementing those provided by the
original proof.

For additional results of Chase and
others on maps from direct product modules
to direct sums, see \cite{SUC1}, \cite{SUC2}, \cite{JO'N},
\cite{BZH}, and papers referred to in the latter two works.
The more recent works obtain results on maps from products
of not necessarily countable families of modules.
It would be interesting to know whether similar results can be
obtained for maps on inverse limits of modules with respect
to not necessarily countable inversely directed systems.

Incidentally, reading Chase~\cite{SUC3} led me to strengthen my
own results by replacing an original descending chain of right ideals
with a downward directed system $\mathcal F,$ and to explicitly state
Theorem~\ref{T.MNB} rather than passing directly to the case where
the $M_i$ are finitely generated, i.e., Corollary~\ref{C.fg}.

In~\cite{BZH}, the right ideals occurring these results are
generalized still further, the operation of multiplying by such an
ideal being replaced
with any subfunctor of the forgetful functor from left
$\!R\!$-modules to abelian groups that commutes with arbitrary direct
products, and it is noted that these include not only functors of
multiplication by finitely generated right ideals, but also functors
obtainable from those by transfinite iteration.

Where above we have examined homomorphisms from direct products and
related constructions to direct sums,~\cite{EK} investigates
homomorphisms from direct products and
related constructions to general modules, and~\cite{EKN}
homomorphisms from general modules to direct sums and related
constructions; though in exchange for this greater generality,
the authors of those papers study a more restricted set of questions.
It is amusing that where, above, I took an argument by successive
approximation and replaced it with an application of the
Baire Category Theorem, the authors of~\cite{EK} take a similar
argument, presented in~\cite{Lady}
in terms of the Baire Category Theorem, and translate it
back into a construction by successive approximation.
{\em Plus \c{c}a change~\dots\,.}

Some similar results with nonabelian groups in place of
modules are proved in \cite{GH}, \cite{KE}, \cite{KE+SSh}.

\section{Generating products of general algebras.}\label{S.gen}

The diagonal argument by which we proved Theorem~\ref{T.fgorunc}
used nothing specific to modules.
Also, as noted at the beginning
of~\S\ref{S.diag}, the focus there on countable products was
for the sake of presentational simplicity.
We shall give the result here in its natural generality.
(However, for consistency of notation, having begun this paper
in the context of module theory, I will not follow the
general algebra convention of using different symbols for
algebras and their underlying sets.)

As noted in the statement of the next theorem, though the hypothesis
there restricts the {\em arities} of the operations, there is
no restriction on the cardinality of the {\em set} of operations.
That cardinality corresponds to the cardinality of the ring $R$ in
Theorem~\ref{T.fgorunc}.
The latter theorem is trivial for $R$ of
cardinality $<2^{\aleph_0},$ and this one is likewise trivial
if there are $<2^\kappa$ operations.

The second conclusion of the theorem, concerning ultraproducts,
was pointed out to me by T.\,Scanlon.
(That conclusion formally subsumes the first conclusion, and of course
implies the corresponding intermediate statements with $U$ replaced by
any filter containing no set of cardinality $<\kappa,$ since any
such filter extends to an ultrafilter with the same property.)
\begin{theorem}\label{T.>*kgens}
Let $\kappa$ be an infinite cardinal, and $T$ a type of
algebras such that all operations of $T$ have arities $<\kappa,$
and if $\kappa$ is singular, all have arities $\leq$ some common
cardinal $<\kappa;$ but where no assumption is made on the cardinality
of the set of those operations.
Let $(M_i)_{i\in\kappa}$ be a $\!\kappa\!$-tuple of algebras
of type $T,$ each of which requires at least $\kappa$ generators.
Then $\prod_{i\in\kappa}M_i$ requires $>\kappa$ generators.

In fact, if $U$ is any ultrafilter on $\kappa$ which
contains no set of cardinality $<\kappa,$ then the ultraproduct
$(\prod_{i\in\kappa}M_i)/U$ requires $>\kappa$ generators.
\end{theorem}\begin{proof}
The restriction on arities of operations of $T$
insures that if an algebra of type $T$
is generated by a set $X,$ then each element of that algebra
belongs to a subalgebra generated by $<\kappa$ elements of $X,$
and in fact that any family of $<\kappa$ elements is contained
in such a subalgebra.
Given this observation, the proof of the first
assertion is exactly analogous to that of Lemma~\ref{L.Pinfg}.

For $X\subseteq \prod M_i$ of cardinality $\leq \kappa,$ $y$ an element
of $\prod M_i$ constructed, as in that proof, to avoid
the subalgebra generated by $X,$ and $z$ any element of the
latter subalgebra, note that $y$ and $z$ will in fact agree
at fewer than $\kappa$ coordinates.
Hence for an ultrafilter $U$ containing no set of cardinality $<\kappa,$
the image of $y$ in $(\prod M_i)/U$ will not equal the image of $z;$ so
$(\prod M_i)/U$ is not generated by the image of such a set $X;$
i.e., it, too, requires $>\kappa$ generators.
\end{proof}

We can now generalize Theorem~\ref{T.fgorunc}.
The value of $\lambda$ of greatest interest in the next result is, of
course, $\aleph_0.$

\begin{theorem}\label{T.bigorsmall}
Let $\lambda\leq\kappa$ be infinite cardinals, with $\lambda$
regular, let $T$ be an algebra type such that all operations of
$T$ have arities $<\lambda,$ and let
$M$ be any algebra of type $T.$
Then the algebra $M^\kappa$ either
requires $<\lambda$ or $>\kappa$ generators.
\end{theorem}\begin{proof}
Let $\mu$ be the least cardinality of a generating set for $M^\kappa.$
If $\mu$ were neither $<\lambda$ nor $>\kappa,$ then noting
that $M^\kappa=(M^\kappa)^\mu,$ we could apply Theorem~\ref{T.>*kgens},
with $\mu$ in place of $\kappa$ and $M^\kappa$ in place of each
$M_i,$ to conclude that $M^\kappa$ requires $>\mu$ generators,
a contradiction.
\end{proof}

Unlike Theorem~\ref{T.>*kgens}, Theorem~\ref{T.bigorsmall} contains
no statement about ultrapowers, so we ask

\begin{question}\label{Q.ultra}
Is the analog of Theorem~\ref{T.bigorsmall} true with $M^\kappa$
replaced by the ultrapower $M^\kappa/U,$ for $U$ any
ultrafilter on $\kappa$ containing no set of cardinality $<\kappa$?
\end{question}

One can get an affirmative answer to this question when $\lambda=\kappa$
for ultrafilters $U$ having a very special property.
For each $\alpha\in\kappa$ let us write
$i_\alpha:\kappa\rightarrow\kappa\times\kappa$ for the
injection $i_\alpha(\beta)=(\alpha,\beta),$ and for any ultrafilter $U$
on $\kappa,$ let
us define an ultrafilter $U^2$ on $\kappa\times\kappa$ by
\begin{xlist}\item\label{x.U^2}
$U^2~=\ \{S\subseteq\kappa\times\kappa~|
\ \{\alpha\in\kappa~|\ i_\alpha ^{-1}(S)\in U\}\in U\}.$
\end{xlist}
Now if $U$ is an ultrafilter on $\kappa$ containing
no set of cardinality $<\kappa,$ and if there exists an injection
$\phi: \kappa\times\kappa\rightarrow\kappa$ such that the
ultrafilter on $\kappa$ gotten by pushing $U^2$ forward
via $\phi$ is $U,$ then one finds that $(M^\kappa/U)^\kappa/U\cong
M^{\kappa\times\kappa}/U^2\cong M^\kappa/U$ (making use of $\phi$
at the last step), and one deduces from Theorem~\ref{T.>*kgens}
that the ultrapower $M^\kappa/U$ must require either $<\kappa$
or $>\kappa$ generators.
However, I do not know whether there are many, or indeed,
any nonprincipal ultrafilters for which such a $\phi$ exists.

Another way one might be able to get a partial or complete affirmative
answer to Question~\ref{Q.ultra} would be to prove some result to the
effect that for fixed $M$ and $\kappa,$ the smallest cardinality of
a set of generators required by an ultrapower $M^\kappa/U$ does not
depend very strongly on the ultrafilter $U.$
If such a result is true (and I have no idea whether one
is), then rather than needing an ultrafilter that is isomorphic
to its own ``square'' in the above argument, it would suffice to have
a family of ultrafilters that ``agree'' in this respect, and
such that the ``product'' in the above sense of some
member of that family with some other was again in the family.

I have no reason to expect positive answers to the next
question; the present wording is just the shortest way of asking for
counterexamples that should exist for part~(i), and may
also for part~(ii), but that I don't know how to find.

\begin{question}\label{Q.inf}
In the final conclusion of Theorem~\ref{T.bigorsmall},\\
\textup{(i)} can one strengthen ``$\!\<<\lambda\!$''
to ``$\!\<<\aleph_0\!$''?\\
\textup{(ii) (}cf.\ Question~\ref{Q.c'uum}\textup{)}
can one strengthen ``$>\kappa\!$'' to ``$\geq 2^\kappa\!$''?
\end{question}

Let us note a ``difficulty'' with Theorem~\ref{T.bigorsmall}
as a generalization of Theorem~\ref{T.fgorunc}: Classes of algebras with
infinitely many operations are not commonly considered in most fields
other than ring- and module-theory.
We will now note one close analog of the module-theoretic situation;
then, in the next section, introduce a class of examples of a quite
different sort, and obtain some results on these.

The nonlinear analog of an $\!R\!$-module is
a $\!G\!$-set where $G$ is a group or monoid, and the next lemma, which
is proved just like Lemma~\ref{L.R=R^*k}, shows
that in that context, the ``$\!<\lambda\!$'' case
of Theorem~\ref{T.bigorsmall} again occurs.
In the category $\mathcal C$ referred to in the statement,
morphisms are assumed to
compose like functions written on the left of their arguments,
i.e., the composite of morphisms $f:X\rightarrow Y$
and $g: Y\rightarrow Z$ is written $gf: X\rightarrow Z.$

\begin{lemma}\label{L.M=M^*k}
Let $\mathcal C$ be a category, $\kappa$ a cardinal, and
$M$ an object of $\mathcal C$ which is either a $\!\kappa\!$-fold
coproduct of copies of itself
\begin{xlist}\item\label{x.M=(+)Mgen}
$M~\cong~\coprod_\kappa M$
\end{xlist}
{\rm (}for instance, if $\kappa$ is infinite and $M$ any object
of the form $\coprod_\kappa N),$ or
a $\!\kappa\!$-fold product of copies of itself
\begin{xlist}\item\label{x.M=M^gen}
$M~\cong~\prod_\kappa M$
\end{xlist}
{\rm (}for instance, if $\kappa$ is infinite and $M$ any object
of the form $\prod_\kappa N).$
In the former case, if we let $G$ be the monoid ${\rm End}(M),$ or in
the latter, if we let $G~=~{\rm End}(M)^{\rm op},$ then
\begin{xlist}\item\label{x.G=G^*k}
$G~\cong~G^\kappa$ as left $\!G\!$-sets.
\end{xlist}
\vspace{-10pt}\qed
\end{lemma}

The situation is strikingly different when $G$ is a group, no matter
how constructed.
If $M$ is any left $\!G\!$-set of more than one element, then for
every proper nonempty subset $S\subseteq\kappa,$ the set $A_S$ of
elements of $M^\kappa$ having exactly two values, one at all
indices in $S$ and the other at all indices in $\kappa-S,$ is a union
of $\!G\!$-orbits, and sets $A_{S_1}$ and $A_{S_2}$ are disjoint unless
$S_1$ and $S_2$ are equal or complements of one another.
So for $\kappa$ infinite, $M^\kappa$ consists of $\geq 2^\kappa$
orbits, and so cannot be generated by fewer elements.

The concept of a $\!G\!$-set, for $G$ a monoid or a group, is
generalized by that of an algebra $M$ of any type on which
$G$ acts by endomorphisms, respectively by automorphisms.
$\!R\!$-modules for $R$ a ring are a particular class of examples;
there may be others for which Theorem~\ref{T.bigorsmall}
would be of interest, and in particular, where the unexpected
phenomenon of generation of $M^\kappa$ by $<\lambda$ elements occurs.

\section{Generation of power algebras over their diagonal
subalgebras.}\label{S.Sym}
A different way to get natural examples of algebras with
operation-sets of arbitrarily large cardinalities is to
start with algebras of arbitrary type $T,$ fix one such
algebra $X$ whose underlying set has large cardinality, and consider
algebras $Y$ given with a homomorphism $i$ of $X$ into them, formally
treating the image of each element of $X$ as a zeroary operation.
The simplest case is that in which $Y=X$ and $i$ is the identity map.
For that case, Theorem~\ref{T.bigorsmall} (with $\lambda$ taken
to be $\aleph_0$ for simplicity) says

\begin{corollary}\label{C.overdiag}
Let $\kappa$ be a cardinal, and $M$ an algebra with
operations all of finite arity.
Then the algebra $M^\kappa$ is either generated by the
diagonal $\Delta(M)$ and finitely many additional elements,
or requires $>\kappa$ additional elements.\qed
\end{corollary}

The next result gives, for any infinite
cardinal $\kappa,$ nontrivial algebras of two
familiar sorts having only finitely many operations, whose
$\!\kappa\!$th direct powers are finitely generated
over their diagonal subalgebras.
In the proof, we continue to write maps on the left of their
arguments and compose them accordingly, though this reverses the
convention in the material from \cite{Sym_Omega:1} to be cited
in the proof of the second statement.

\begin{theorem}\label{T.Sym}
Let $\kappa$ be an infinite cardinal and $\Omega$ a set
of cardinality $\geq\kappa.$\\
{\rm (i)}  If $M$ is the monoid of all maps $\Omega\rightarrow\Omega,$
then $M^\kappa$ is generated over $\Delta(M)$ by two elements.\\
{\rm (ii)}  If $S$ is the group of all permutations of $\Omega,$
then $S^\kappa$ is generated over $\Delta(S)$ by one element.
\end{theorem}\begin{proof}
(i):  Since $|\Omega|=\kappa\cdot|\Omega|,$ we can write $\Omega$ as the
union of $\kappa$ disjoint sets of the same cardinality
as $\Omega,$ $\Omega=\bigcup_{i\in\kappa}\Sigma_i.$
For each $i\in\kappa,$ let $a_i\in M$ be a one-to-one map
with image $\Sigma_i,$ and $b_i\in M$ a left inverse to $a_i.$
I claim $M^\kappa$ is generated by $\Delta(M)$ and the
two elements $a=(a_i)_{i\in\kappa}$ and $b=(b_i)_{i\in\kappa}.$
For given any $x=(x_i)_{i\in\kappa}\in M^\kappa,$ let us ``encode''
$x$ in a single element $x'\in M,$ defined to act on each
subset $\Sigma_i$ by $a_i\,x_i\,b_i.$
Then we see that $b\,\Delta(x')\,a=(b_i\,x'\,a_i)_{i\in\omega}=
(b_i\,a_i\,x_i\,b_i\,a_i)_{i\in\omega}=(x_i)_{i\in\omega}=x,$
as required.

(ii):  Again the idea will be to encode elements of $S^\kappa$ in single
elements of $S.$
The trouble is that our structure no longer contains elements
$a_i$ and $b_i$ giving bijections between $\Omega$
and the subsets of $\Omega$ on which we will encode the
components of our $\!\kappa\!$-tuple.
It will, however, contain permutations carrying each of these subsets
bijectively to a common set $\Sigma_1\subset\Omega.$
But note that if we take an element whose actions on
these subsets ``encode'' the coordinates of a member of $S^\kappa,$
and conjugate it by each of these permutations to get permutations
having those same actions on
$\Sigma_1,$ the information encoded in its other components
will not disappear; it will move to other parts of $\Omega,$ where it
will constitute ``garbage'' that we must get rid of.
We will do this using commutator operations, in which those
parts of our map are commuted with identity maps.
Finally, we will call on a result from \cite{Sym_Omega:1} to go from
permutations of $\Sigma_1$ to permutations of $\Omega.$

So let $\Omega=\Sigma_1\cup\Sigma_2,$ where $\Sigma_1$ and $\Sigma_2$
are disjoint sets of the same cardinality as $\Omega,$ and
let $S_1$ be the subgroup of $S$ consisting of elements which
fix all members of $\Sigma_2$ and act arbitrarily on $\Sigma_1.$
We shall first show that $S_1^\kappa$ is contained in the subgroup of
$S^\kappa$ generated over $\Delta(S)$ by a single
element $f=(f_i)_{i\in\kappa}.$

Let us identify $\Sigma_2$ with $\Sigma_1\times\kappa\times
\{0,1\},$ as we may since $|\Omega|\geq\kappa.$
For each $i\in\kappa,$ we define $f_i\in S$ so that it cyclically
permutes the three sets $\Sigma_1,$ $\Sigma_1\times\{i\}\times\{0\}$
and $\Sigma_1\times\{i\}\times\{1\}$
while for $j\in\kappa-\{i\}$ it interchanges the two sets
$\Sigma_1\times\{j\}\times\{0\}$ and $\Sigma_1\times\{j\}\times\{1\}.$
Precisely, we let
\begin{xlist}\item\label{x.fdef}
$f_i:\ \alpha\mapsto(\alpha,i,0)\mapsto(\alpha,i,1)\mapsto\alpha,$~
$(\alpha,j,0)\leftrightarrow(\alpha,j,1),$\\
\indent for $\alpha\in\Sigma_1,\ j\neq i$ in $\kappa.$
\end{xlist}

Now in the group of permutations of an infinite
set, every element is a commutator by \cite{OO} or
\cite[Lemma~14]{Sym_Omega:1}, so given $x=(x_i)_{i\in\omega}$ in
$S_1^\kappa,$ let us, for each $i\in\kappa,$ regard
$x_i$ as a permutation of $\Sigma_1$ and write it as the commutator of
permutations $u_i,\,v_i$ of that set.
Let us now define elements $y,z\in S$ which act trivially
on $\Sigma_1$ and on $\Sigma_1\times\kappa\times\{1\},$
while on each set $\Sigma_1\times\{i\}\times\{0\},$
$y$ behaves like $u_i$ and $z$ like $v_i,$ i.e.,
$y(\alpha,i,0)=(u_i(\alpha),i,0),\ z(\alpha,i,0)=(v_i(\alpha),i,0).$
I claim that $x$ is the commutator of $f^{-1}\Delta(y)f$ and
$f^2\Delta(z)f^{-2}.$
Indeed, looking at the $\!i\!$th component of $f^{-1}\Delta(y)f$
for any $i,$ we see that it acts trivially except on $\Sigma_1$ and on
the sets $\Sigma_1\times\{j\}\times\{1\}$ with $j\neq i,$
while $f^2\Delta(z)f^{-2}$ is trivial except on $\Sigma_1$ and on
the sets $\Sigma_1\times\{j\}\times\{0\}$ with $j\neq i.$
Hence their commutator is trivial except on $\Sigma_1,$ where it
behaves as the commutator of $u_i$ and $v_i,$ i.e., as $x_i.$
So the subgroup of $S^\kappa$ generated by $\Delta(S)\cup\{f\}$
contains every $x\in S_1^\kappa,$ as claimed.

For the remainder of the proof we put aside the description of
$\Omega$ as $\Sigma_1\cup(\Sigma_1\times\kappa\times
\{0,1\})$ used to get this fact, and take a decomposition of
a simpler sort, keeping $\Sigma_1$ as above, but
choosing a set $\Sigma_3$ such that
$\Sigma_1\cup\Sigma_3=\Omega,$ and such that
$\Sigma_1\cap\Sigma_3,$ $\Sigma_1-\Sigma_3$ and $\Sigma_3-\Sigma_1$
all have cardinality $|\Omega|.$
From the proof of~\cite[Lemma~2]{Sym_Omega:1} one finds that in this
situation, if we denote by $S_3$ the group of permutations of $\Omega$
that fix all elements not in $\Sigma_3,$ then every
permutation of $\Omega$ can be written either as a member of the
product-set $S_1\,S_3\,S_1$ or of the product-set
$S_3\,S_1\,S_3;$ hence we have $S=S_1\,S_3\,S_1\,S_3.$
(See last three paragraphs proof of~\cite[Lemma~2]{Sym_Omega:1}.
If one doesn't want to read between the lines of
that proof, one can use the statement of that lemma,
which, for subsets $U$ and $V$ of $S$ satisfying certain weaker
conditions than being equal to the above subgroups, shows --
after adjusting to left-action notation -- that every element of $S$
belongs either to $V (V\,U)^4$ or to $U (U\,V)^4.$
Hence if $U$ and $V$ contain $1,$ every element will belong
to $(U\,V)^6,$ hence in our present
situation, to $(S_1\,S_3)^6,$ which one may use
in place of $S_1\,S_3\,S_1\,S_3$ in the reasoning of the next sentence.)
Now taking $t\in S$ which interchanges $\Sigma_1$ and $\Sigma_3,$
we get $S_3=t^{-1}S_1\,t,$
so $S=S_1\,(t^{-1}S_1\,t)\,S_1\,(t^{-1}S_1\,t),$
hence $S^\kappa = S_1^\kappa\,\Delta(t^{-1})S_1^\kappa\,
\Delta(t)\,S_1^\kappa\,\Delta(t^{-1})S_1^\kappa\,\Delta(t)$.
Since the subgroup of $S^\kappa$ generated by $\Delta(S)\cup\{f\}$
contains $S_1^\kappa,$ the above equation shows that it is
all of $S^\kappa.$
\end{proof}

The result from \cite{Sym_Omega:1} called on in the last
paragraph of the above proof was used there in proving
two other properties of infinite symmetric
groups, essentially~(\ref{x.unccof}) and~(\ref{x.bddgen}) below.
One may ask whether there is a direct relation between
the conclusion of the above lemma and those properties.
An implication in one direction is proved, under simplifying
restrictions on the value of $\kappa$ and the algebra-type, in

\begin{theorem}\label{T.fg/=>}
Let $S$ be an algebra with finitely many primitive operations,
all of finite arity, which satisfies
\begin{xlist}\item\label{x.fg/*D}
The countable power algebra $S\<^\omega$ is finitely generated
over the diagonal subalgebra $\Delta(S).$
\end{xlist}
Then $S$ also satisfies both of
\begin{xlist}\item\label{x.unccof}
$S$ cannot be written as the union of a countable chain of
proper subalgebras.
\end{xlist}
\begin{xlist}\item\label{x.bddgen}
For every subset $X\subseteq S$ which generates $S,$ there exists a
positive integer $n$ such that all elements of $S$ can be
represented by words of length $\leq n$ in the elements of $X.$
\end{xlist}
\end{theorem}\begin{proof}
It is not hard to show that the conjunction of~(\ref{x.unccof})
and~(\ref{x.bddgen}) is equivalent to the statement
\begin{xlist}\item\label{x.unc+bdd}
Whenever $X_0\subseteq X_1\subseteq\dots$ is an $\!\omega\!$-indexed
chain of subsets of $S$ with $\bigcup_{i\in\omega} X_i=S,$
such that for each primitive operation
$f$ of $S$ and each $i\in\omega,$
one has $f(X_i,\dots,X_i)\subseteq X_{i+1},$
then $X_i=S$ for some $i\in\omega.$
\end{xlist}
So it suffices to show that~(\ref{x.fg/*D}) implies~(\ref{x.unc+bdd}).
(This reduction only needs the fact that~(\ref{x.unc+bdd}) implies
both~(\ref{x.unccof}) and~(\ref{x.bddgen}).
These implications are seen by taking for $X_i,$ in the first case,
the $\!i\!$th member of an $\!\omega\!$-indexed ascending chain of
subalgebras, and in the second, the set of elements of $S$ that can be
represented by words of depth $\leq i$ in the elements of $X.)$
We will prove this in contrapositive form, showing that if
$(X_i)_{i\in\omega}$ is a family which satisfies the
hypotheses of~(\ref{x.unc+bdd}) but not the conclusion, and
$Y\subseteq S\<^\omega$ is finite, then the subalgebra
of $S\<^\omega$ generated by $\Delta(S)\cup Y$ must be proper.

Indeed, define the {\em rank} of an element $s\in S$ to be the
least $r$ such that $s\in X_r.$
Then our assumption that the conclusion
of~(\ref{x.unc+bdd}) fails says that the rank function is unbounded,
so for each $i\in\omega$ we may take for $x_i$ an element of rank
at least $i+{\rm max}_{y\in Y}\,{\rm rank}(y_i)$ (where
each $y\in Y\subseteq S\<^\omega$ is written $(y_i)_{i\in\omega}).$
We claim that the element $x=(x_i)_{i\in\omega}$ does not lie in the
subalgebra of $S\<^\omega$ generated by $\Delta(S)\cup Y.$

For if it did, it would lie in the subalgebra
generated by $\Delta(Z)\cup Y$ for some finite
subset $Z\subseteq S,$ and be represented by a word of
some depth $d$ in these elements; thus for each
$i,$ $x_i$ would be expressed as a word of depth $d$ in the
elements of $Z$ and the $\!i\!$th components of the elements of $Y.$
We now get a contradiction on taking any
$i>d+{\rm max}_{z\in Z}\,{\rm rank}(z),$ since then by choice of $x_i,$
\begin{xlist}\item[]
${\rm rank}(x_i)\ \geq\ i+{\rm max}_{y\in Y}\,{\rm rank}(y_i)\\
\indent
>\ d+{\rm max}_{z\in Z}\,{\rm rank}(z)+
{\rm max}_{y\in Y}\,{\rm rank}(y_i),$
\end{xlist}
contradicting the existence of such a depth-$\!d\!$ expression
for $x_i.$
\end{proof}

Conditions~(\ref{x.unccof}) and, more recently,~(\ref{x.bddgen}) have
been proved for several sorts of groups that arise in ways similar to
infinite symmetric groups (see \cite[\S3]{Sym_Omega:1} for references),
and also for the endomorphism ring of a direct sum
or product of infinitely many copies of any nontrivial module \cite{ZM}.
It seems likely that some or all of those proofs can be
modified to establish~{\rm (\ref{x.fg/*D})} in these cases as well.
The endomorphism-ring case in fact follows easily from
Lemma~\ref{L.R=R^*k} above and the observation that for any ring $R,$
\begin{xlist}\item\label{mod>bi>rg}
$R\<^\omega$ fin.gen.\ as left $\!R\!$-module $\implies$
$R\<^\omega$ fin.gen.\ as $\!(R,R)\!$-bimodule\\
\indent
$\implies$ $R\<^\omega$ fin.gen.\ over $\Delta(R)$ as ring,
\end{xlist}
where the last step uses the fact that in
$R\<^\omega,$ left or right module multiplication by an element of $R$
is equivalent to ring multiplication by an element of $\Delta(R).$
If in this argument we replace rings and modules with monoids and
$\!M\!$-sets, and the application of Lemma~\ref{L.R=R^*k}
with an application of Lemma~\ref{L.M=M^*k}, with
$\mathcal C=\mathcal{S}et,$ we get Theorem~\ref{T.Sym}(i)
with ``two elements'' strengthened to ``one element''.
(I stated Theorem~\ref{T.Sym}(i) as I did because the proof
given seemed the best way to lead up to the more
difficult proof of part~(ii).
In proving directly the ``one element'' statement, one can take
that one element to be what we called $a$ in the earlier proof.
Given $x=(x_i),$ if we define $x''\in M$ to have restriction to
each $\Sigma_i$ given by $x_i b_i,$ we see that $\Delta(x'')a=x.)$

A nontrivial {\em finite} algebra clearly satisfies~(\ref{x.unccof})
and~(\ref{x.bddgen}), but will not satisfy~(\ref{x.fg/*D}),
since $\Delta(S)$ is finite while $S^\omega$ is uncountable;
so the converse of Theorem~\ref{T.fg/=>} is false.
On the other hand (still considering only algebras
with finitely many operations, all of finite arities), we see that
a finitely generated infinite algebra
cannot satisfy~(\ref{x.bddgen}), and a countably generated but not
finitely generated algebra cannot satisfy~(\ref{x.unccof}),
so any infinite algebra satisfying both~(\ref{x.unccof})
and~(\ref{x.bddgen}) must be uncountable, raising the hope
that such algebras might also satisfy~(\ref{x.fg/*D}).
However, this is not so,
for it is shown in \cite{YdC} and in \cite[Corollaire~18]{AK} (both
generalizing a result in \cite{SK+JT})
that an infinite direct power of a finite perfect group
satisfies~(\ref{x.unccof}) and~(\ref{x.bddgen}); but such a group cannot
satisfy~(\ref{x.fg/*D}), since it admits a homomorphism
onto a nontrivial finite group, and~(\ref{x.fg/*D}) clearly
carries over to homomorphic images.

Since the proof of Theorem~\ref{T.fg/=>} constructs an element $x$
which disagrees at all but finitely many coordinates with
every element of the subalgebra generated by $\Delta(S)\cup Y,$
one can in fact say (using again the idea of the last
sentence of Theorem~\ref{T.>*kgens})
that for any nonprincipal ultrafilter $U,$ the image of such
an element $x$ in $S^\kappa/U$ fails to lie in the image
of the subalgebra generated by $\Delta(S)\cup Y.$
Thus, for $U$ a nonprincipal ultrafilter on $\omega,$
one can insert in Theorem~\ref{T.fg/=>} the intermediate condition
\begin{xlist}\item\label{x.fg//U}
The ultrapower $S\<^\omega/U$ is finitely generated
over the image of the diagonal subalgebra $\Delta(S),$
\end{xlist}
and strengthen that theorem
to say that (\ref{x.fg/*D})$\Rightarrow$(\ref{x.fg//U})$\Rightarrow$%
(\ref{x.unccof})$\wedge$(\ref{x.bddgen}).

Unlike~(\ref{x.fg/*D}), but like~(\ref{x.unccof}) and~(\ref{x.bddgen}),
condition~(\ref{x.fg//U}) is clearly satisfied by finite algebras.
For groups, it is also preserved by the operation of pairwise direct
product, since that operation respects
the operation ``subgroup generated by''.
Hence the direct product of the permutation
group on an infinite set with any nontrivial finite group is an
example of an infinite algebra that separates~(\ref{x.fg/*D})
from~(\ref{x.fg//U}).
The example mentioned earlier of an infinite direct power of a
finite perfect group turns out to separate~(\ref{x.fg//U}) from
(\ref{x.unccof})$\wedge$(\ref{x.bddgen}), in view of the next result.

\begin{proposition}\label{P.xPfin}
If an algebra $S$ satisfies~{\rm (\ref{x.fg//U})} for some
nonprincipal ultrafilter $U$ on $\omega,$ then there is
a {\rm (}finite{\rm )} upper bound on the cardinalities of
finite homomorphic images of $S.$
Equivalently, $S$ has a least congruence $C$ such that $S/C$ is finite.
\end{proposition}\begin{proof}
Let $S$ be an algebra which does not have a least congruence $C$
making $S/C$ finite.
We shall show that $S$ does not satisfy~(\ref{x.fg//U}).

Denote by $\mathcal F$ the class of all congruences $C$ on $S$ with
finite quotient $S/C,$ and let
$C_{\rm res.fin.}=\bigcap_{C\in\mathcal F} C.$
Since~(\ref{x.fg//U}) is preserved under passing to homomorphic
images, we may divide out by $C_{\rm res.fin.}$ and
assume $S$ residually finite, but still infinite.

I claim now that we can find a chain of congruences
$C_0\supset C_1\supset\dots\in\mathcal F$ and a sequence of
elements $u_0,\,u_1,\,\dots\in S$ such that
$(u_i,\,u_{i+1})\in C_i-C_{i+1}$ for all $i.$
Indeed, let $C_0\in\mathcal F$ be arbitrary.
Since $S/C_0$ is finite, some congruence class $c_0$ with respect to
$C_0$ is infinite.
Choose a congruence $C_1\subset C_0$ such that $c_0$
decomposes into more than one congruence class under $C_1.$
At least one of these must be infinite; choose such an infinite
class $c_1\subset c_0,$ and repeat the process.
For each $i,$ choose $u_i\in c_i-c_{i+1};$
then the relations $(u_i,\,u_{i+1})\in C_i-C_{i+1}$ clearly hold.
Writing $p_i$ for the canonical map $S\rightarrow S/C_i,$
these relations say that for all~$i,$
\begin{xlist}\item\label{x.cong+not}
$p_{i+1}(u_i)~\neq\ p_{i+1}(u_{i+1})~=\ p_{i+1}(u_{i+2})~
=\ p_{i+1}(u_{i+3})~=\ \dots\,.$
\end{xlist}

To show that~(\ref{x.fg//U}) fails, let
$Y\subseteq S\<^\omega$ be any finite subset;
we must show that the image in $S\<^\omega/U$ of the subalgebra
of $S\<^\omega$ generated by $\Delta(S)\cup Y$ is a proper subalgebra.

Since $S/C_0$ and $Y$ are finite, and $U$ is an
ultrafilter, there exists $R_0\in U$ such that for each
$y=(y_j)_{j\in\omega}\in Y,$
the $\!\omega\!$-tuple $(p_0(y_j))_{j\in\omega}\in (S/C_0)^\omega$
becomes constant when $j$ is restricted to $R_0\subseteq\omega.$
Likewise, we can find $R_1\subseteq R_0$ in $U$ such that
for each $y\in Y,$ $(p_1(y_j))_{j\in\omega}$ becomes constant when
restricted to $j\in R_1,$ and so on,
getting a chain $R_0\supseteq R_1\supseteq\dots$ of members of
$U$ such that for each $i\in\omega$ and each $y\in Y,$
the elements $p_i(y_j)\in S/C_i$ are the same for all $j\in R_i.$
Also, since $U$ is nonprincipal, we can along the way make sure, for
each $i,$ that $i\notin R_i,$ so that $\bigcap_i R_i=\emptyset.$

Now the condition that the $\!j\!$- and $\!j'\!$-components of
all elements of $Y$ have equal images under $p_i$ for certain
$i,$ $j,$ and $j'$ carries over to elements of
the subalgebra generated by $\Delta(S)\cup Y;$ hence
\begin{xlist}\item\label{x.ij=ij'}
For all $z=(z_j)$ in the subalgebra of $S\<^\omega$ generated by
$\Delta(S)\cup Y,$ all $i\in\omega,$ and all $j,j'\in R_i,$
we have $p_i(z_j)=p_i(z_{j'}).$
\end{xlist}
We shall now construct an element $x\in S\<^\omega$ which does not
agree modulo $U$ with any element $z$ satisfying~(\ref{x.ij=ij'}),
in other words, such that for any such $z,$ and any $R\in U,$
\begin{xlist}\item\label{x.notinR}
$\{j\in\omega~|\ x_j=z_j\}\ \not\subseteq\ R.$
\end{xlist}

For each $j\in\omega,$ let $i(j)$ be the
greatest integer such that $j\in R_{i(j)}$ if $j\in R_0,$ and $-1$
otherwise; and let $x_j$ be $u_{i(j)-1}$ if $i(j)\geq 1,$
arbitrary otherwise.
To prove~(\ref{x.notinR}) given $R\in U,$ take any $j\in R_1\cap R,$
and then any $j'\in R_{i(j)+1}\cap R,$ noting that $1\leq i(j)<i(j').$
By~(\ref{x.ij=ij'}), with $i=i(j),$
we have $p_{i(j)}(z_j)=p_{i(j)}(z_{j'}).$
On the other hand, $x_j$ and $x_{j'}$ are
$u_{i(j)-1}$ and $u_{i(j')-1}$ respectively,
so~(\ref{x.cong+not}) with $i=i(j)-1$ shows
that $p_{i(j)}(x_j)\neq p_{i(j)}(x_{j'}).$
Thus, the $\!j\!$th and $\!j'\!$th components of $x$ cannot
{\it both} coincide with the corresponding components of $z.$
Since $j$ and $j'$ both belong to $R,$ this proves~(\ref{x.notinR}),
as required.
\end{proof}

Since we do not know the answer to Question~\ref{Q.ultra},
even for $\lambda=\omega=\aleph_0,$ we in particular
do not know whether~{\rm(\ref{x.fg//U})} is equivalent to
\begin{xlist}\item\label{x.cg//U}
The ultrapower $S\<^\omega/U$ is {\em countably} generated
over the image of the diagonal subalgebra $\Delta(S),$
\end{xlist}
so it is worth noting that in both the implication
(\ref{x.fg//U})$\Rightarrow$(\ref{x.unccof})$\wedge$(\ref{x.bddgen})
and the above proposition, the hypothesis~(\ref{x.fg//U})
can be weakened to~(\ref{x.cg//U}), with some extra work.
(Key changes:  If $Y$ is countable, write it as the union
of a chain $Y_0\subseteq Y_1\subseteq\dots$ of finite subsets.
In the construction of the element $x$ as in the proof
of Theorem~\ref{T.fg/=>}, for each $i$ let $x_i$ have rank
$\geq i+{\rm max}_{y\in Y_i}\,{\rm rank}(y_i),$ using the
finite set $Y_i$ rather than all of $Y.$
Likewise, in the proof of Proposition~\ref{P.xPfin}, construct the
sets $R_i$ so that the elements $p_i(y_j)\in S/C_i$ $(j\in R_i)$ are
merely the same for each $y\in Y_i.)$

It is natural to ask

\begin{question}\label{Q.U,U'}
For a given algebra $S,$ if~{\rm(\ref{x.fg//U})} holds for
some nonprincipal ultrafilter $U,$ will it hold for all
nonprincipal ultrafilters?
If not, what implications can be obtained among these conditions for
different $U$?
\end{question}

Khelif~\cite[D\'efinition~5]{AK} introduces a condition having an
interesting similarity to~(\ref{x.fg/*D}).
Given an algebra-type $T$ and a natural number $n,$ let $W(n,T)$ denote
the set of all $\!n\!$-ary words
in the operations of $T.$
Note that if $S$ is an algebra of type $T,$
then an element of $W(n,T)^\omega$ induces an $\!n\!$-ary
operation on $S\<^\omega.$
Khelif's definition (reformulated to bring out the parallelism
with~(\ref{x.fg/*D})) says that an algebra $S$ ``has property P*\,'' if
there exists a natural number $n$ and an element
$w=(w_i)_{i\in\omega}\in W(n,T)^\omega$ such that every element
$x\in S\<^\omega$ is the value of $w$ at some $\!n\!$-tuple of
elements of $\Delta(S).$
In~\cite[D\'efinition~4]{AK} he defines a weaker condition P, saying
that there exists a natural
number $n,$ and a function $\eta$ from the
natural numbers to the natural numbers, such that for every
$x\in S\<^\omega$ there exists $w\in W(n,T)^\omega$ whose $\!i\!$th
component is a word length $\leq\eta(i)$ for each $i,$ such that,
again, $x$ is the value of $w$ at an $\!n\!$-tuple of
elements of $\Delta(S);$ and he shows that P (and hence
also P*) implies (\ref{x.unccof})$\wedge$(\ref{x.bddgen}).

Since there are only finitely many words of a given length, and
there is a bound on the lengths of any finite set of words, Khelif's
property P is equivalent to the statement that for some $n$ one can
associate to every natural number $i$ a finite subset
of $W_i\subseteq W(T,n)$ such that every element
$x\in S\<^\omega$ is the value of some $w\in\prod_i W_i$
at some $\!n\!$-tuple of elements of $\Delta(S).$
This suggests a parallel generalization of~(\ref{x.fg/*D}):
\begin{xlist}\item\label{x.pfg/*D}
There exists a sequence $(X_i)_{i\in\omega}$ of finite
subsets $X_i\subseteq S$ such that $S\<^\omega$ is generated
by $\Delta(S)\,\cup\,\prod_i\,X_i.$
\end{xlist}
And indeed, one sees that the proof of Theorem~\ref{T.fg/=>}
(and its strengthening by the insertion of conditions~(\ref{x.fg//U})
and~(\ref{x.cg//U}))
works equally well under this weaker hypothesis, since the construction
there of the sequence $(x_i)$ called on the finiteness of the
generating set $Y$ only via the finiteness of the
subsets $\{y_i~|~y\in Y\}$ $(i\in\omega).$

One can also combine features of~(\ref{x.fg/*D}) or~(\ref{x.pfg/*D}) and
Khelif's conditions, allowing, for instance, finitely
many elements of $W(T,n)^\omega$ and finitely many
elements of $S\<^\omega$ to be used together in generating
$S\<^\omega$ over $\Delta(S\<^\omega).$
This plethora of conditions leads one to wonder whether there is
some small number of ``core'' conditions to which most or all of
those we have mentioned are related in simple ways; e.g., by adding
conditions such as ``$\!S\!$ has no finite homomorphic images'',
adding cases such as ``all finite algebras'', and/or varying
some natural parameters in the conditions.

It is noted in \cite[paragraph following Question~8]{Sym_Omega:1}
that no non-finitely-generated abelian group satisfies
either~(\ref{x.unccof}) or~(\ref{x.bddgen}), so
Theorem~\ref{T.fg/=>} shows (in two ways) that
no nontrivial abelian group satisfies~(\ref{x.fg/*D}).
One may ask generally

\begin{question}\label{Q.vars}
What can be said about varieties of algebras containing
non-finitely-generated algebras $S$ that satisfy
{\rm (\ref{x.fg/*D})}?
{\rm (\ref{x.unccof})}?
{\rm (\ref{x.bddgen})}?
{\rm (\ref{x.fg//U})}?
{\rm (\ref{x.cg//U})}?
Khelif's P*?
P?
Further variants of these?
\end{question}

Another condition, weaker than both~(\ref{x.bddgen}) and
Khelif's property P, and of a more elementary nature,
which might be of interest in the study
of these conditions, is the statement that for every generating
set $X$ for $S,$ there exist an integer $n$ such that every
element of $S$ belongs to a subalgebra generated by $\leq n$
elements of $X.$
This is satisfied by the abelian group $Z_{p^\infty},$
with $n=1$ independent of $X.$

Turning in a different direction, observe that if an algebra $S$
satisfies~(\ref{x.fg/*D}), then so will $S\<^\omega.$
For if we write $S\<^\omega=S',$ and identify $(S')^\omega$ with
$S^{\omega\times\omega}\cong S\<^\omega,$ then the diagonal
subalgebra of $(S')^\omega$ contains the diagonal subalgebra of
$S^{\omega\times\omega},$ over which~(\ref{x.fg/*D}) shows
that it is finitely generated.
In particular, $S\<^\omega$ satisfies~(\ref{x.bddgen}); so given
a finite set $Y$ such that $\Delta(S)\cup Y$ generates $S\<^\omega$
as in~(\ref{x.fg/*D}), there will be a bound on the lengths of
words needed to so express every element of $S\<^\omega.$
It would be interesting to know whether one can in general
get all elements of $S\<^\omega$ using a
single word, with specified positions in each of which a
specified element of $Y$ occurs, while arbitrary elements of $\Delta(S)$
are put in the other positions, as the proof of Theorem~\ref{T.Sym}
shows that one can in the case of the monoid or group of all
maps or invertible maps of an infinite set into itself.

To an algebra $S$ satisfying~(\ref{x.fg/*D}), we can associate several
natural-number-valued invariants witnessing that condition:
The least cardinality of
a set $Y$ such that $S\<^\omega$ is generated by $\Delta(S)\cup Y;$
the least $n$ such that for some finite $Y\in S\<^\omega$ every element
of $S\<^\omega$ can be expressed using words of length (or depth)
$\leq n$ in terms of elements of $\Delta(S)\cup Y,$ etc..
Note that if we have a sequence $S_0,\ S_1,\,\dots$ of algebras of the
same type which each satisfy~(\ref{x.fg/*D}), but for which
the sequence of values of
such an invariant is unbounded, then the product algebra
$\prod_{i\in\omega} S_i$ cannot satisfy~(\ref{x.fg/*D}),
in contrast to the observation of the preceding
paragraph on a direct power of one such $S.$

A confession about Theorem~\ref{T.fg/=>}:  The hypothesis that $S$ had
only finitely many primitive operations was not used in the proof.
I made that assumption because it would be needed in places in our
subsequent discussion, and because without it, the
conclusion~(\ref{x.bddgen}) of that theorem is {\em weaker}
than optimal.
The reader will not find it hard to verify that if $S$ has
{\em countably} many primitive operations, then
for any generating set $X,$ all elements of $S$ may be
obtained not merely using words of finite length as in~(\ref{x.bddgen}),
but using such words in finitely many of these operations,
and more generally, that if $S$ has an arbitrarily large set of
primitive operations, and this set is
partitioned in any way into countably many subsets, then $S$
can be obtained as in~(\ref{x.bddgen}) using the operations
in the union of finitely many of these subsets.
(Key idea:  In~(\ref{x.unc+bdd}), replace the
relation $f(X_i,\dots,X_i)\subseteq X_{i+1}$
by $f(X_i,\dots,X_i)\subseteq X_{i+a(f)},$ where $a$ is a function
from the set of primitive operations to the positive integers giving the
partition into countably many subsets, and in place of ``depth'' in the
proof, use a ``weighted depth'' that takes account of this function.
Incidentally, this result is analogous to~(\ref{x.unccof}), since the
latter can be translated
as saying that whenever a generating set for $S$ is partitioned into
countably many subsets, the union of some finite number of these
subsets is still a generating set.)
For an example of an algebra $S$ satisfying~(\ref{x.fg/*D})
with uncountably many primitive operations, such that no
{\em finite} set of these suffices, let $S=\omega,$ and let the set of
primitive operations be $\omega^\omega,$ i.e., the set of
{\em all} unary operations on $S.$
Then $S\<^\omega$ is generated under these
operations by a single element, the identity
function, but one cannot obtain all of $S$ from the generating
set $\{0\}$ using expressions of bounded length in finitely many
operations.\vspace{6pt}

Some final miscellaneous remarks:

An easy class of examples of algebras that satisfy~(\ref{x.bddgen})
but not~(\ref{x.unccof}) is that of infinite algebras in which
the value of each operation is always one of its arguments; for
instance, an infinite chain regarded as a lattice, or
an infinite set with no operations.
These examples also show that if we generalize
\cite[Question~8]{Sym_Omega:1},
which asked whether a countably infinite group can
satisfy~(\ref{x.bddgen}), to algebras of arbitrary type,
the answer is affirmative.

Although~(\ref{x.fg/*D}) and~(\ref{x.unccof}) are preserved under
taking homomorphic images, and although the same is true
of~(\ref{x.bddgen}) in any algebra whose structure
includes a structure of group, it is not true
of~(\ref{x.bddgen}) in general, for if $S'$ is a homomorphic
image of $S,$ the inverse image of a generating set for $S'$
may not be a generating set for $S.$
For example, let $S'$ be the semilattice of finite subsets
of $\omega$ under union, and $S$ the semilattice obtained by
hanging from each $x\in S',$ a new element
\begin{xlist}\item\label{x.x0<x}
$x_0\ <x\ ,$
\end{xlist}
with no new
order-relations other than the consequences of~(\ref{x.x0<x}),
so that, in particular, the new elements are pairwise incomparable.
The retraction $S\rightarrow S'$ taking each $x_0$ to $x$
is easily seen to be a homomorphism,
so $S'$ is a homomorphic image of $S.$
Now any generating set for $S$ must contain the set of join-irreducible
elements, $\{\emptyset\}\cup\{x_0~|\ x\in S'\},$ and every element
of $S$ is a join of at most two members of this set, so $S$
satisfies~(\ref{x.bddgen}); but~(\ref{x.bddgen}) fails for $S',$
as shown by the generating set consisting of $\emptyset$ and
all singletons.

Condition~(\ref{x.bddgen}) has been called by some authors ``the Bergman
property'', based on its introduction in~\cite{Sym_Omega:1}.
Another possible name for an algebra with this property is an
{\em impatient} algebra, since it
says ``If you're going to generate me, you have to
do it in finitely many steps -- I can't wait forever!''

Though~(\ref{x.fg/*D}) is, for simplicity's sake, formulated
above only for the $\!\omega\!$-fold direct power, the corresponding
conditions on general powers $S^\kappa$ are, of course, of interest.
(It implies the $\!\kappa\!$-analog of~(\ref{x.unccof}),
called ``cofinality $>\kappa$''.)

And, of course, the diagonal subalgebra $\Delta(S)$
of a power $S^\kappa$ was merely the most obvious example of the
context introduced in the first paragraph of this section,
which would be worth considering in greater generality.

Another kind of algebra structure involving an arbitrarily large
set of operations, all of finite arity, is a set on which
a Boolean ring $B$ ``acts'', in the sense of \cite{bool_act}.
So Theorem~\ref{T.bigorsmall} also applies to these structures,
and it might be of interest to examine that case,
possibly combining the Boolean operations with others.

\end{document}